\documentclass[final]{siamltex}

\usepackage{graphicx} 
\usepackage{amsfonts,amsmath,amssymb}
\usepackage{booktabs,array,multirow}
\usepackage{cancel}
\usepackage[nocompress]{cite}
\usepackage{color}
\usepackage[font=small,labelfont=bf]{caption}
\usepackage{url}
\DeclareMathOperator{\spn}{span}
\newcommand{\argmin}[1]{\underset{#1}{\mathrm{argmin}}}

\usepackage{algorithm}
\usepackage{algpseudocode}

\makeatletter
\let\c@algorithm\c@figure
\makeatother

\title{fastGMRES: Robust randomized GMRES \\ 
via flexible GMRES} 

\author{Stefan G\"{u}ttel\thanks{Department of Mathematics, The University of Manchester, Oxford Road, Manchester, M13\,9PL, United Kingdom, \texttt{stefan.guettel@manchester.ac.uk}}       \and
John W.~Pearson\thanks{School of Mathematics and Maxwell Institute for Mathematical Sciences, The University of Edinburgh, James Clerk Maxwell Building, The King's Buildings, Edinburgh, EH9\,3FD, United Kingdom, \texttt{j.pearson@ed.ac.uk}}}

\begin{document}

\maketitle

\begin{abstract}
We propose a practical randomized algorithm for solving large sparse systems of linear equations. The algorithm, termed fastGMRES, uses flexible GMRES as an outer wrapper for sketched GMRES together with a practical residual bound to terminate the inner iteration. We demonstrate that fastGMRES requires very little parameter tuning, while still being efficient and robust in the sense of generating non-increasing residual norms. We observe that the inner sketched GMRES method can be run without any inner product computations and vector--vector updates on high-dimensional vectors, yielding significant speedup for many problems.
\end{abstract}



\pagestyle{myheadings}
\thispagestyle{plain}
\markboth{S. G\"{U}TTEL AND J. W. PEARSON}{FASTGMRES: ROBUST RANDOMIZED GMRES}

\section{Introduction}

Randomization has become a key ingredient in some of the most efficient numerical linear algebra solvers, in particular, for low-rank matrix approximation and highly overdetermined least squares problems \cite{sarlos2006improved,rokhlin2008fast,WoolfeLibertyRokhlinTygert2008,woodruff2014sketching,martinsson2020randomized}. For these problems, randomization can lead to impressive performance improvements while retaining desirable accuracy and stability properties. The resulting algorithms are hence well suited for implementation in state-of-the-art high-performance computing libraries~\cite{avron2010blendenpik,murray2023randomized}.

More recently, work in the field has been dedicated to speeding up Krylov solvers for  systems of linear equations $Ax=b$, where $A\in\mathbb{R}^{n\times n}$ is a large sparse matrix and $b,x\in\mathbb{R}^n$ \cite{BalabanovGrigori22,NakatsukasaTropp21,BGS23,jang2025randomized,bucci2025randomized}. It is much harder to extract significant performance gains through randomization for this type of problems. Intuitively, this can be explained by the fact that there is no redundancy that can easily be exploited: the solution vector $x$ of a general nonsymmetric system of linear equations can be arbitrarily sensitive to perturbations in any of $A$ or~$b$'s components, hence may depend sensitively on all of them. As a consequence,  instead of following the ``sketch-the-problem'' paradigm as in the case of overdetermined least squares problems, randomized Krylov solvers for $Ax=b$ try to exploit randomization on the level of the algorithm (``sketch-the-algorithm''). 
One of the dominating computational kernels in common Krylov solvers such as GMRES~\cite{SaadSchultz1986}, GMRES-DR~\cite{Morgan2002}, or GCRO-DR~\cite{parks2006recycling}, is the orthogonalization of a Krylov basis in each cycle. 
Randomization has been used in \cite{NakatsukasaTropp21,BalabanovGrigori22,BGS23} to reduce these expensive orthogonalization costs.

Let us consider a linear system $A x = b$ with initial guess $x_0$ and associated residual $r_0 := b - A x_0$. 
We will assume that the columns of $V_{m}\in\mathbb{R}^{n\times m}$ span some search space. 
We further assume that we have a \emph{sketching operator} $S\in\mathbb{R}^{s\times n}$ with $m < s\ll n$ which acts as an approximate isometry for the Euclidean norm $\|\cdot\|$. More precisely, given a positive integer $m$ and some $\varepsilon\in [0,1)$, let $S$ be such that for all vectors~$v\in\spn(V_m)$,
\begin{equation}
(1-\varepsilon) \| v \|^2 \leq \| S v\|^2 \leq (1+\varepsilon) \|v\|^2.
\label{eq:sketch}
\end{equation}
The mapping $S$ is called an \emph{$\varepsilon$-subspace embedding} for $\spn(V_m)$; see, e.g.,~\cite{sarlos2006improved,woodruff2014sketching,martinsson2020randomized}. 
Condition~\eqref{eq:sketch} can equivalently be stated with the Euclidean inner product~\cite[Cor.~4]{sarlos2006improved}: for all $u,v \in \spn(V_m)$,  
\begin{equation}{\label{eq:sketch_innerproduct}}
\langle u, v \rangle - \varepsilon \| u\|\cdot \|v\|
            \leq \langle Su, Sv \rangle 
            \leq \langle u, v \rangle + \varepsilon \| u\|\cdot \|v\|.\nonumber
\end{equation}
In practice, $S$ is not explicitly available, but we can draw it at random to achieve~\eqref{eq:sketch} with high probability.

In the case of \emph{sketched GMRES} (sGMRES)~\cite{NakatsukasaTropp21}, the columns of $V_{m}$ form a (not necessarily orthogonal) basis of the Krylov space $$\mathcal{K}_m(A,r_0) := \spn\{r_0,Ar_0, \ldots,A^{m-1}r_0\}.$$ The algorithm computes a correction $V_{m} y_{m}$ such that the norm of the associated sketched residual $S r_{m}$ is minimal in the Euclidean norm, leading to the sGMRES approximation
\begin{equation}\label{eq:sGMRES_problem}
  x_{m} = x_{0} + V_{m} y_{m}, \quad    
  y_{m} = \argmin{y \in \mathbb{R}^m}
        \|
             S r_{0} - S A V_{m} y 
        \|.
\end{equation}
The size of the sketched least squares problem in~\eqref{eq:sGMRES_problem} depends only on $s$ and~$m$, independent of the problem size~$n$. It can thus be solved cheaply without requiring the computation of an orthonormal basis~$V_{m}$ as in standard GMRES. 

The central goal in this work is to devise a method that is reliable and robust, is competitive with state-of-the-art Krylov methods for a wide range of problems, and (unlike a number of such techniques) requires very little user input in the form of parameter tuning. With this in mind, we advocate for the use of flexible GMRES (FGMRES; see \cite{saad1993flexible}) as a stabilizing method for sGMRES. Here, ``stabilizing'' should be understood synonymously with ``robustifying'' the method, in the sense that the combination of FGMRES and sGMRES becomes more reliable while retaining efficiency. In particular, the outer FGMRES iteration guarantees non-increasing residual norms in exact arithmetic and reduces stagnation issues often encountered with restarted \mbox{sGMRES}. 

The outline of this paper is as follows. In Section~\ref{sec:fgmres} we briefly review the FGMRES method and introduce a flexible version of the full orthogonalization method, FFOM. This is then used to derive a practical upper bound on the FGMRES residual, and provide new insights on the convergence of the method. Section~\ref{sec:stab} then combines these developments into the proposed FGMRES-sGMRES routine, termed fastGMRES. Numerical experiments are reported in Section~\ref{sec:numex}, before we conclude and discuss possible future work in Section~\ref{sec:concl}.

\section{FGMRES and FFOM}\label{sec:fgmres}

For ease of reference, we state in Algorithm~\ref{alg:fgmres} the FGMRES algorithm from~\cite{saad1993flexible}, as well as a slight modification we refer to as \emph{flexible full orthogonalization method} (FFOM). FFOM can be seen as a natural extension of the FOM method \cite[Sec.~6.4]{saad2011numerical} to variable preconditioning. 

We note that in this section we use $e_*$ to denote the standard basis vectors, the dimension of which is to be inferred from the context. Throughout, $x_m/r_m$ denote the $m$-th FGMRES iterate/residual, $\widetilde x_m/\widetilde r_m$ are the $m$-th FFOM iterate/residual, respectively, and $\widehat r_m$ will denote the residual of the inner FGMRES/FFOM preconditioner.

\begin{algorithm}[h!]
\caption{Pseudocode for the FGMRES and FFOM algorithms.}\label{alg:fgmres}
\begin{algorithmic}
\Require Matrix $A\in\mathbb{R}^{n\times n}$, vector $b\in\mathbb{R}^n$, initial guess $x_0\in\mathbb{R}^n$, $m\geq 1$
\State $r_0 := b - Ax_0$, $\beta:=\|r_0\|$, $v_1:=r_0/\beta$
\For{$j=1,\ldots,m$}
    \State $z_j := M_j^{-1} v_j$\quad 
(Section~\ref{sec:stab}: Solve $A z_j \approx v_j$ with sGMRES, zero initial guess)
    \State $w := A z_j$
    \For{$i = 1,\dots,j$}
        \State $h_{i,j} := v_i^* w$; \ $w :=  w - h_{i,j}v_i$
    \EndFor
    \State $h_{j+1,j} := \| w\|$; \ $v_{j+1} := w / h_{j+1,j}$
\EndFor
\State Define $Z_{m} := [ z_1,z_2,\ldots,z_{m} ]$
\State \textbf{FGMRES:} Compute $x_m:= x_0 + Z_m y_m$, where $y_m$ minimizes $\|\beta e_1 - \underline{H_m} y \|$. That is, the residual norm $\|r_m\| = \| b - A x_m \|$ is minimized over all $x_m= x_0 + Z_m y_m$.
\State \textbf{FFOM:} Compute $\widetilde x_m:= x_0 + Z_m y_m$ with $y_m = H_m^{-1}\beta e_1$. We define $\widetilde r_m := b - A \widetilde x_m$.
\end{algorithmic}
\end{algorithm}

Both algorithms produce basis matrices $V_m\in\mathbb{R}^{n\times m}$ and $Z_m\in\mathbb{R}^{n\times m}$, with the vectors in $Z_m$ being the preconditioned counterparts of the vectors in $V_m$; that is, $z_j := M_j^{-1} v_j$. 
While FGMRES is designed to minimize the residual norm $\|r_m\| = \| b - A x_m \|$ over all $x_m = x_0 + Z_m y_m$, FFOM instead enforces that the residual $\widetilde r_m$ be orthogonal to the span of $V_m$. 
For computing their approximants, both algorithms exploit an Arnoldi-like decomposition of the form $A Z_m = V_{m+1} \underline{H_m}$, where the $(m+1)\times m$ upper Hessenberg matrix $\underline{H_m}$ collects the orthogonalization coefficients $h_{i,j}$. The upper $m\times m$ part of $\underline{H_m}$ will be denoted by $H_m$. 
Trivially, we have
\[
\|r_m\| = \| b - A x_m \| \leq \| b - A \widetilde x_m \| = \| \widetilde r_m \|
\]
for the FGMRES residual $r_m$ and the FFOM residual $\widetilde r_m$, and so in this sense there appears to be no obvious reason to choose FFOM over FGMRES. However, FFOM admits a concise characterization of the residual $\widetilde r_m$ as
\begin{align*}
b - A \widetilde x_{m} 
&=  b - A x_0 - A Z_{m} H_m^{-1} \beta e_1 \\
&=  r_0 - V_{m+1} \underline{H_m} H_{m}^{-1} \beta e_1 \\
&=  r_0 - (V_{m} H_{m} + v_{m+1} h_{m+1,m}e_{m}^T) H_{m}^{-1} \beta e_1 \\
&= - v_{m+1} h_{m+1,m} \underbrace{e_{m}^T H_{m}^{-1} \beta e_1}_{=:\gamma_m},
\end{align*}
and so
\begin{equation}
\| \widetilde r_m\| = \| b - A \widetilde x_{m} \| = h_{m+1,m}\cdot |\gamma_{m}|.
\label{eq:ffomres}
\end{equation}

\subsection{A practical bound on the FGMRES residual} \label{sec:bnd}

The operator $M_j^{-1}$ in Algorithm~\ref{alg:fgmres} is an (inverse) preconditioning operation that can change at each iteration. Frequently, $M_j \approx A$ and we would like to understand how the approximation quality of the preconditioner affects the residual decrease of FGMRES. 

Let us start with the first iteration $j=1$ and assume that $z_1 = A^{-1} (v_1 - \widehat r_1)$. That is, the preconditioner solves the linear system $A z_1 = v_1 - \widehat r_1$, where $\widehat r_1$ is a (small) residual. 
As $Z_1 = [z_1]$, we have
\[
r_1 = b - A x_1 = b - A(x_0 + z_1 \gamma_1)
= r_0 - \gamma_1 (v_1 - \widehat r_1) 
= \beta v_1 - \gamma_1 (v_1 - \widehat r_1),
\]
with the coefficient $\gamma_1$ chosen by FGMRES so that $\|b - A x_1\|$ is minimized. 
If we instead choose  $\gamma_1 = \beta = \| r_0 \|$, we get the upper bound
\[
\|r_1\| = \| b - A x_1 \| \leq \beta \| \widehat r_1 \|,
\]
that is, the first FGMRES iteration reduces the initial residual norm $\beta = \| r_0\|$ by at least the factor~$\| \widehat r_1 \|$.

Now assume that at iteration $j=m>1$ we have $z_m = A^{-1} (v_m - \widehat r_m)$. 
Using the relation $AZ_m = V_{m+1} \underline{H_m}$, we have
\begin{align*}
   r_m =  b - A x_m 
    &= b - A (x_0 + Z_m y_m)  \\
    &= r_0 - A Z_{m-1} \begin{bmatrix}
    \gamma_1\\
    \gamma_2\\
    \vdots\\
    \gamma_{m-1}
    \end{bmatrix}
    - (v_m - \widehat r_m) \gamma_m \\
    &= 
    r_0 - V_m \underline{H_{m-1}} \begin{bmatrix}
    \gamma_1\\
    \gamma_2\\
    \vdots\\
    \gamma_{m-1}
    \end{bmatrix} 
    - (v_m - \widehat r_m) \gamma_m,
\end{align*}
with the coefficients $\gamma_j$  chosen by FGMRES so that $\| b - A x_m \|$ is minimized. 

An alternative, convenient choice of coefficients $\gamma_j$ is when 
\[
\underline{H_{m-1}} \begin{bmatrix}
    \gamma_1\\
    \gamma_2\\
    \vdots\\
    \gamma_{m-1}
    \end{bmatrix} 
    = 
\begin{bmatrix}
    \beta\\
    0\\
    \vdots\\
    0\\
   - \gamma_m
    \end{bmatrix}.
\]
With this choice of coefficients, 
\[
r_m = b - A x_m = \gamma_{m} \cdot \widehat r_m.
\]
Assuming $H_{m-1}$ is invertible, the solution is 
\[
\begin{bmatrix}
    \gamma_1\\
    \gamma_2\\
    \vdots\\
    \gamma_{m-1}
\end{bmatrix}  
=    
H_{m-1}^{-1} 
\begin{bmatrix}
    \beta\\
    0\\
    \vdots\\
    0
    \end{bmatrix}, \quad \gamma_m = -h_{m,m-1}\cdot \gamma_{m-1}.
\]
Since FGMRES minimizes the residual over all possible coefficient choices, we conclude that
\begin{equation}
\|r_m \| = \| b - Ax_m \| \leq h_{m,m-1}\cdot |\gamma_{m-1} | \cdot \|\widehat r_m\|.
\label{eq:fgmres_bnd}
\end{equation}
Note that $h_{m,m-1}$ and $\gamma_{m-1}$ are available \emph{before} $z_m$ is computed, and hence we can use this bound to predict the FGMRES residual decrease based on the residual norm $\|\widehat r_m\|$ of the preconditioner. If the preconditioner at iteration $j=m$ is exact ($\widehat r_m=0$), the FGMRES iterate $x_m$ will also be exact. This has already been known~\cite[Sec.~2.2]{saad1993flexible}.

We further observe that $h_{m,m-1}\cdot |\gamma_{m-1}|$ is in fact the residual norm of the $(m-1)$-th FFOM approximant 
\[
\widetilde x_{m-1} = x_0 + Z_{m-1} H_{m-1}^{-1} \beta e_1;
\]
see \eqref{eq:ffomres}. We can therefore state \eqref{eq:fgmres_bnd} more concisely as $\|r_m\| \leq \|\widetilde r_{m-1}\| \cdot \|\widehat r_m\|$, or perhaps more verbatim,
\[
\boxed{\big\| r_m^{\text{FGMRES}} \big\| \leq 
\big\| r_{m-1}^{\text{FFOM}} \big\|\cdot 
\big\| r_{m}^{\text{PRECOND}} \big\|.}
\]
(This bound is also essentially contained in the proof of~\cite[Lem.~3]{vuik1995new}.) 
Therefore, as long as the flexible FOM and flexible GMRES approximants at iteration $m-1$ (that is, $x_{m-1}$ and $\widetilde x_{m-1}$) have similar residual norms, we can say that the FGMRES residual at iteration $m$ is reduced approximately by the factor $\| \widehat r_m \|$, i.e., the residual achieved by the preconditioner.

\begin{figure}
\centering\includegraphics[width=0.9\linewidth]{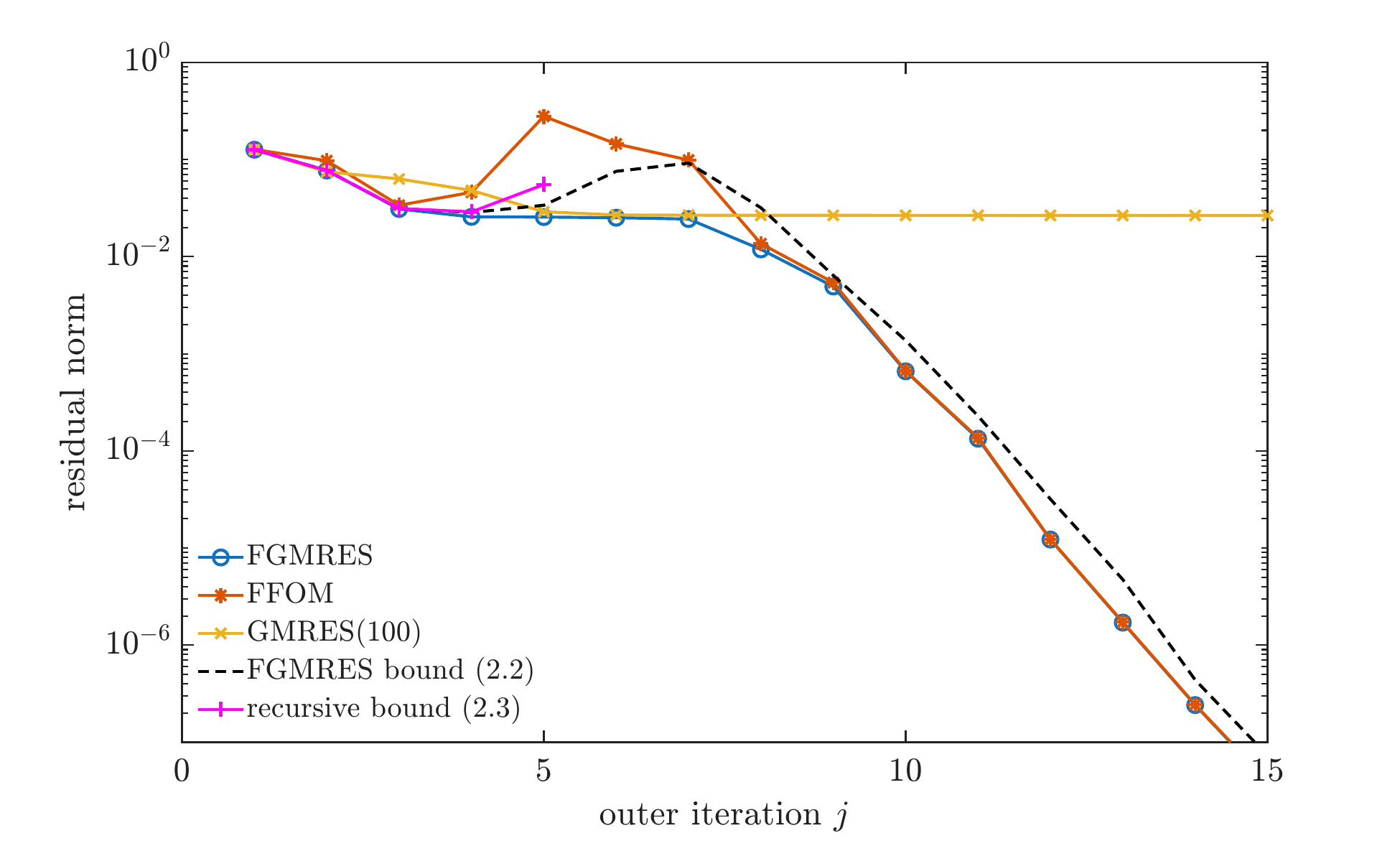}
    \caption{Illustration of FGMRES and FFOM convergence on an artificial matrix of size $n=1,000$. The inner (inverse) preconditioning operation $M_j^{-1}$ corresponds to performing 100 GMRES iterations. We also show the convergence of restarted GMRES(100) as a standalone solver. The dashed black line shows our upper bound~\eqref{eq:fgmres_bnd} on the FGMRES residual. At each iteration $j$, this bound can be evaluated cheaply by using only the residual of the $(j-1)$-st FFOM approximant and the residual of the inner preconditioner $M_j$. The recursive bound \eqref{eq:lauri} is shown as well (it breaks down after 5 iterations).}
    \label{fig:fgmres_bnd}
\end{figure}

Figure~\ref{fig:fgmres_bnd} demonstrates the FGMRES residual bound \eqref{eq:fgmres_bnd} for an artificial example. Here, the matrix $A$ of size $n=1,000$ has been generated in MATLAB via \texttt{A = randn(n)+30*eye(n)} and \texttt{b = randn(n,1)}. The inner solver $M_j^{-1}$ corresponds to performing 100 GMRES iterations. This is of comparable leading-order arithmetic cost to performing restarted GMRES with a restart length of 100, apart from the additional outer FGMRES operations. However, as can be seen in Figure~\ref{fig:fgmres_bnd}, restarted GMRES stagnates and fails to converge for this problem.
Both FGMRES and FFOM, on the other hand, converge steadily and the FGMRES residual is closely tracked by our residual upper bound.

\subsection{FGMRES versus restarting}

As an alternative to preconditioning, we can simply think of $M_j^{-1}$ as an iterative solver that is restarted using the residual--error relation $Ae = r$. This gives rise to the pseudocode in Algorithm~\ref{alg:restprecon}, which can be considered as applying ``iterative improvement'' with a variable approximation ($M_j$) of $A$. Restarted GMRES falls into that class of methods.

\begin{algorithm}[h!]
\caption{Restarting with a variable preconditioner $M_j$.}\label{alg:restprecon}
\begin{algorithmic}
\Require $A\in\mathbb{R}^{n\times n}$, right-hand side $b\in\mathbb{R}^n$, initial guess $x_0\in\mathbb{R}^n$, $m\geq 1$
\State $r_0 := b - Ax_0$
\For{$j=1,\ldots,m$}
    \State $z_j := M_j^{-1} r_{j-1}$
    \State $x_j := x_{j-1} + z_j$
    \State $r_j := b - A x_j$
\EndFor
\end{algorithmic}
\end{algorithm}

Let $\widehat r_j := r_{j-1} - A z_j$ be the residual of the iterative solver in Algorithm~\ref{alg:restprecon}. Then 
\[
r_j = b - A x_j = b - A x_{j-1} - A z_j
    = \widehat r_j,
\]
that is, the residual $r_j$ for the problem $Ax=b$ is precisely equal to the residual~$\widehat r_j$ of the iterative solver. If the iterative solver is GMRES, then~$r_j$ will be minimized locally (that is, within each restart cycle).  FGMRES, on the other hand, performs a residual minimization step over all previously computed vectors in~$Z_m$.

We highlight that the above does \emph{not} imply that FGMRES with $k$~iterations of GMRES as the inner solver necessarily achieves a smaller residual norm than GMRES($k$) at each cycle. It is possible to construct examples where standalone GMRES($k$) converges faster than FGMRES with GMRES($k$) as the preconditioner.

\subsection{Remarks on convergence}

Generally, the analysis of FGMRES is complicated because this method extracts its approximants from a subspace spanned by $Z_m$ and that space is no longer a standard Krylov subspace. As a consequence, results on the convergence of GMRES derived using polynomial approximation theory do not carry over to the flexible setting. However, in the case that the inner preconditioner is a polynomial Krylov method (like GMRES or sGMRES), the search space spanned by $Z_m$ is a subspace of a larger Krylov space. 
This is the setting analyzed in 
\cite{simoncini2002flexible}. It is proven that for any flexible Arnoldi-type method which uses a constant order-$m$ Krylov method as the preconditioner, the space spanned by the union of all Krylov
vectors generated by the preconditioners (denoted $\mathcal{B}_k$ in that paper) keeps growing in dimension. 
That analysis applies in particular to FGMRES-GMRES($m$). 
However, that result does not imply FGMRES convergence nor finite termination after $n$ iterations.

Another discussion of FGMRES can be found in \cite{simoncini2003theory}. Therein, the viewpoint is to interpret FGMRES as an inexact Krylov subspace method. One starts with a non-variable preconditioned system $AM^{-1} (Mx) = b$, but then performs inexact products with variable $A M_j^{-1}$ that approximate $AM^{-1}$. Informally, if $M_j$ and $M$ are sufficiently close at every iteration, then we expect to get essentially the same convergence as for GMRES on the preconditioned system $AM^{-1}$ plus an error term, say, $O(\epsilon)$. Then, if exact GMRES for $AM^{-1} (Mx) = b$ converges geometrically at rate $\rho$, in the inexact case we would expect a residual decrease like
$$\|r_j\| \le C \rho^j \|r_0\| + O(\epsilon).$$
However, this interpretation of FGMRES as an inexact Krylov method is somewhat cumbersome: one would need to quantify how closely the preconditioners $M_j$ relate to $M$ and make sure that the difference is not too large. Also, the $O(\epsilon)$ term suggests that FGMRES may never be able to reduce the residual to zero, which is not usually the case. It therefore seems unlikely that inexact Krylov theory can be used to obtain realistic convergence bounds for FGMRES.

Recent work in \cite{GN26} has studied the convergence of FGMRES under the assumption that the residuals $\widehat r_j$ produced by the preconditioner are bounded in norm, using \eqref{eq:fgmres_bnd} as a starting point. It is shown that, if all $\|\widehat r_j\|\leq \mu \leq 1/2$, then FGMRES is guaranteed to converge with a geometric rate between $\mu$ and $\sqrt{\mu}$. The work also presented a recursive upper bound on the FGMRES residuals $r_m$ given by
\begin{equation}\label{eq:lauri}
\|r_m \| \leq \delta_m \| r_{m-1} \|, \quad \delta_m = \frac{\|\widehat r_m\|}{\sqrt{1-\delta_{m-1}^2}}, \quad \delta_0 = 0.
\end{equation}
Conveniently, this bound does not require the FFOM residuals, but on the other hand it is only valid for all $m$ as long as $\max_{j\leq m} \delta_j <1$. This condition may be difficult to satisfy without strict residual control. An example of this is shown in  Figure~\ref{fig:fgmres_bnd}, where the recursive bound breaks down from iteration~6 (because $\delta_5\approx 1.9$ and $\delta_6$ becomes imaginary).

In summary, \eqref{eq:fgmres_bnd} provides a particularly practical bound on FGMRES convergence for our purposes, and we will base our algorithmic developments on it.

\begin{figure}
\centering\includegraphics[width=0.52\linewidth]{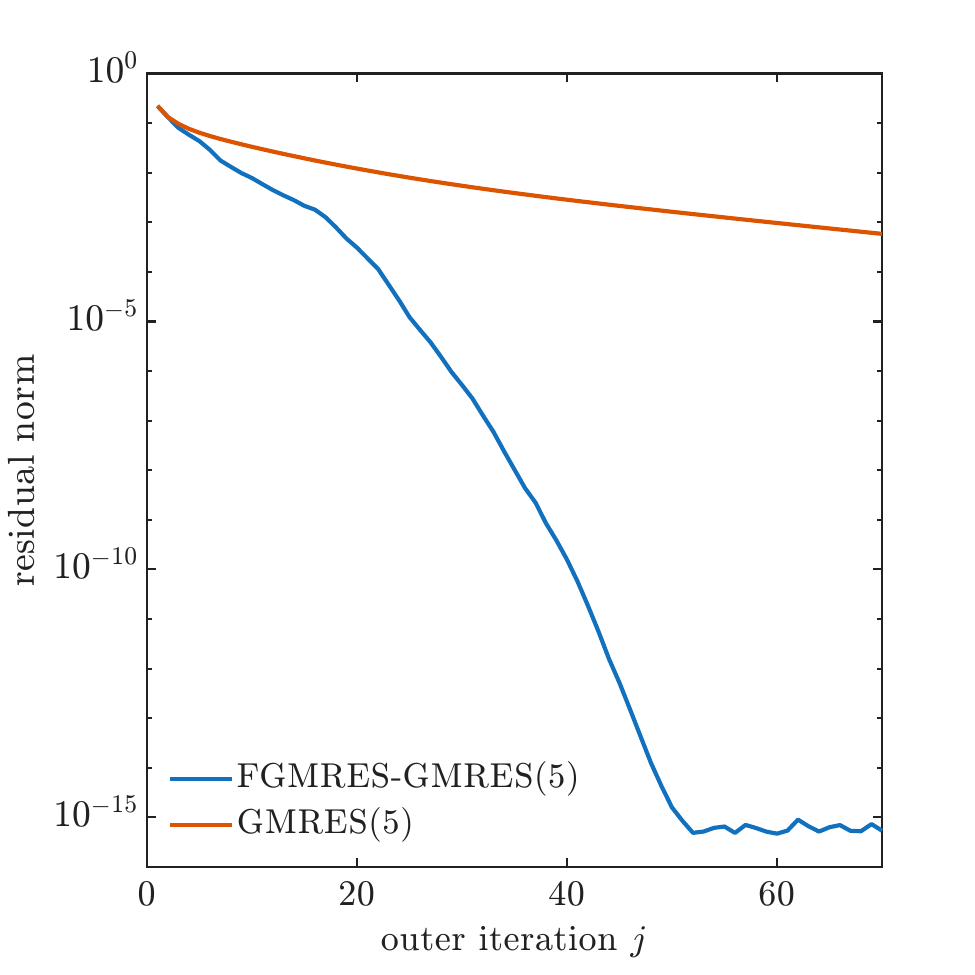}\hspace*{-2mm}\includegraphics[width=0.52\linewidth]{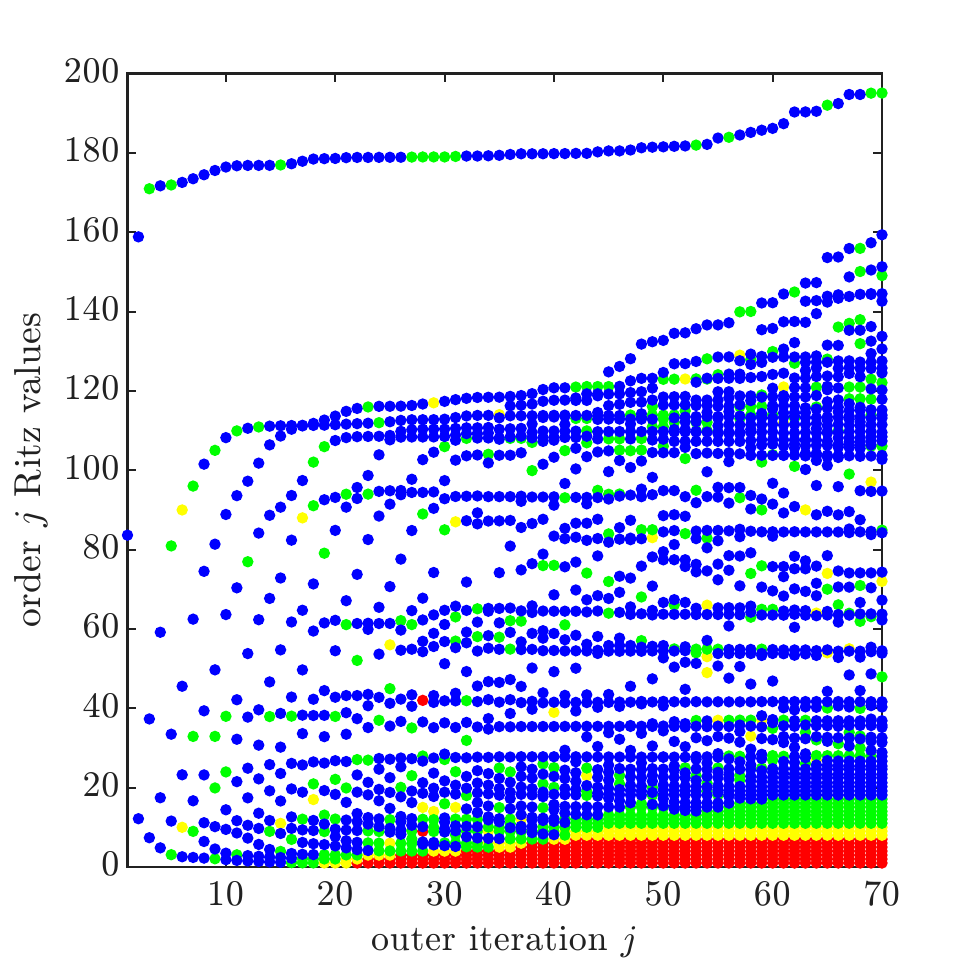}
    \caption{Illustrating spectral deflation effects in FGMRES on a diagonal matrix with eigenvalues $1,2,\ldots,1000$. The preconditioner is GMRES($5$). Left: FGMRES-GMRES($5$) clearly exhibits superlinear convergence in contrast to GMRES($5$). Right: Order $j$ Ritz values associated with the FGMRES search space spanned by $Z_j$ ($j=1,\ldots,70$). The colors indicate distance to a closest eigenvalue of $A$ in decreasing order: blue ($\textrm{distance}\geq 0.1$), green ($0.01 \leq  \textrm{distance} < 0.1$), yellow ($0.001 \leq  \textrm{distance} < 0.01$), and red ($\textrm{distance} < 0.001$).}
    \label{fig:deflation}
\end{figure}

Finally, let us comment on another property of FGMRES that sets it apart from restarted GMRES: potential for increasing spectral deflation with each outer iteration~$j$. If the preconditioners $M_j$ are good approximations to $A$, we expect the vectors $z_j = M_j^{-1} v_j$ to have large components in the direction of eigenvectors associated with small eigenvalues. That is, the search space $Z_m$ from which the FGMRES approximants~$x_m$ are extracted may contain improved approximations of these eigenvectors, leading to convergence acceleration of FGMRES. This is illustrated in Figure~\ref{fig:deflation}. Here, we run FGMRES on a diagonal matrix with eigenvalues $1,2,\ldots,1000$, using GMRES($5$) as the preconditioner. The vector $b$ was chosen at random using MATLAB's \texttt{randn} function. Figure~\ref{fig:deflation} (left) shows the FGMRES residuals, clearly exhibiting superlinear convergence as $j$ increases. By contrast, restarted GMRES($5$) converges only linearly. Figure~\ref{fig:deflation} plots, for each outer FGMRES iteration $j=1,2,\ldots,70$, the order~$j$ Ritz values associated with the FGMRES search space $\text{range}(Z_j)$, that is, the eigenvalues of the matrix $Z_j^\dagger A Z_j$. We observe that for relatively small $j$, some of these Ritz values start converging to the eigenvalues of $A$ closest to zero. This behavior is very similar to convergence of rational Ritz values when the shift parameter of the rational Krylov space is fixed at zero, and the superlinear convergence effects observed when approximating more general matrix functions~$f(A)b$; see, e.g., the plots in \cite{BGV10,BG12}. It is perhaps possible to relate FGMRES to an inexact rational Krylov iteration and thereby prove superlinear convergence in specific settings. We leave this for future work as it is beyond the scope of this paper.

\section{Stabilizing randomized GMRES}\label{sec:stab}

Randomized sketched GMRES may suffer from instabilities or stagnation, for example, in the approach by \cite{NakatsukasaTropp21} where short recurrences (e.g., truncated Arnoldi) were used to construct a possibly highly ill-conditioned Krylov basis. (See also \cite{simoncini2005effect,GS23b} for discussions of nonorthogonal Krylov bases and their effect on convergence.) While \cite[Sec.~3.4]{NakatsukasaTropp21} recommended restarting the sGMRES iteration when the basis condition number exceeds a threshold (like $10^{15}$ in double precision), experiments in \cite{GS23b} revealed that the basis condition number alone may not be a good indicator for the convergence behavior of sGMRES. Rather, it appears that the full set of singular values is important, or at least the numerical rank of the generated basis. Instability issues with too ill-conditioned Krylov bases were also the motivation to develop  randomized GMRES variants with deflation and restarting; see, e.g., \cite{BGS23}. 

Another line of work in \cite{jang2025randomized} proposed randomized variants of GMRES with deflated restarting, as well as a restarted FGMRES method with deflation. In these variants, the generated Krylov bases $V_m$ are kept nearly orthogonal by orthogonalizing their sketches~$S V_m$. This effectively controls the condition number of the Krylov basis itself and hence instabilities are usually not observed. On the downside, generating a near-orthogonal Krylov basis involves all previously generated basis vectors and so the arithmetic cost of this approach still grows quadratically with the number of generated basis vectors. Note that  \cite[Alg.~3]{jang2025randomized} incorporated deflation by effectively restarting the FGMRES method after~$m$ steps, extracting a number of approximate right singular vectors to be projected out of the operator~$A$ in the subsequent cycle. This is different from the GMRES-SDR method of \cite{BGS23}, where a (non-flexible) GMRES method itself was restarted with a sketched variant of the Krylov--Schur algorithm~\cite{Stewart01}. 

Here we propose instead to simply wrap the sGMRES method within an outer FGMRES loop. The global minimization of the outer FGMRES loop then guarantees non-increasing residual norms in exact arithmetic, even if the inner sGMRES method stagnates or becomes unstable. 
The proposed combination of methods is concisely described by Algorithm~\ref{alg:fgmres}, with $M_j^{-1}$ corresponding to $k$ iterations of sGMRES~\cite{NakatsukasaTropp21}. More precisely, at each outer FGMRES iteration~$j$, we compute $z_j$ as an approximate solution of $Az = v_j$:
\begin{itemize}
\item Compute a (non-orthogonal) basis $B_k$ of the Krylov space
$$\mathcal{K}_k(A,v_j)=\mathrm{span}\{ v_j, A v_j, \ldots, A^{k-1} v_j \}.$$
\item Define $z_j$ as the minimizer of $\| S (v_j - A z)\|$ over all $z \in \mathrm{range}(B_k)$, i.e.,
\[
z_j := B_k y_k, \quad y_k := (SAB_k)^\dagger (Sv_j).
\]
\end{itemize}

The number $k$ of sGMRES iterations will be chosen as large as possible, limited by the following three conditions:
\begin{itemize}
\item[\textbf{C1:}] $k\leq k_{\max}$, where $k_{\max}$ is chosen rather large, such as $k_{\max}=500$. This is to limit the required memory for storing the Krylov basis matrix $B_k$, and to allow for an informed choice of the sketching parameter $s$, for example, setting $s := 2 k_{\max}$ as a practical heuristic.
\item[\textbf{C2:}] $\mathrm{cond}(SAB_k)\leq \texttt{threshold}$, with $\texttt{threshold}=10^{15}$ in all our experiments and $\mathrm{cond}(\cdot)$ denoting the condition number of a matrix. Although this condition alone may not be enough to guarantee that sGMRES converges robustly, noting that $\mathrm{cond}(SAB_k)$ is not basis-invariant as it depends on the scaling of the columns of $B_k$, we find the outer FGMRES stabilization will take care of any instabilities so that we find this to be an acceptable criterion. Note that the condition number can be computed cheaply if a QR factorization $Q_k R_k = SAB_k$ is available, and the same factorization can be used to solve the least-squares problem for $y_k$. In our code, we use a Gram--Schmidt procedure with reorthogonalization to update $Q_k R_k$ from one iteration to the next.
\item[\textbf{C3:}] We terminate when $\| r_j\| \leq \texttt{tol}$, where $r_j = b - A x_j$ is the residual of the $j$-th FGMRES iteration. We can use the bound \eqref{eq:fgmres_bnd} to stop the inner sGMRES solver based on the (estimated) sGMRES residual at iteration~$k$, and hence predict an outer residual norm below $\texttt{tol}$.
\end{itemize}

\smallskip

Our justification for \textbf{C2} is that we want to spend as much effort as possible at earlier iterations. That is because the cost of the preconditioner sGMRES grows essentially linearly with the number of iterations~(see also the numerical tests in \cite{NakatsukasaTropp21}), while the expected residual decrease is approximately proportional to $\| r_m^\text{PRECOND}\|$ as per our bound. The only cost that grows quadratically is that of the outer FGMRES loop and hence we want to keep its iteration count as small as possible. We therefore believe that this strategy is close to optimal in reducing the overall cost.

\section{Numerical experiments}\label{sec:numex}

Below we test the proposed FGMRES-sGMRES combination on a number of challenging problems. 
All experiments are run on MATLAB Online with an Intel(R) Xeon(R)
Platinum 8375C CPU @ 2.90GHz, using MATLAB R2025a. 
The codes are available at
\smallskip
\begin{center}
\url{https://github.com/nla-group/fastGMRES}
\end{center}

\subsection{Problem \texttt{vas\_stokes\_1M}} 
We solve $A x = b$, where $A$ is the \texttt{vas\_stokes\_1M} matrix of size $n=1,090,664$  from the SuiteSparse Matrix Collection~\cite{davis2011university}. The right-hand side vector is generated via $\texttt{b = randn(n,1)}$, and we use (zero-fill) ILU left-preconditioning. This problem is known to be hard for restarted GMRES, requiring a restart length of (approximately) $m\geq 400$ for convergence. 

The sketching operator is the Clarkson--Woodruff sketch $S\in \{-1,0,1\}^{s\times n}$, where each column has a single element $\pm 1$ (sign with equal probability) in a random row position; see \cite{woodruff2014sketching}. This sparse sketching operator can be applied to a vector in $O(n)$ operations. The sketching dimension is fixed at $s = 2k_{\max}$, where $k_{\max}=500$ is the maximum number of inner sGMRES iterations. The sGMRES bases $B_k$ are generated using the truncated Arnoldi process; see, e.g., \cite[Alg.~6.6]{Saad2003}, projecting each new Krylov basis vector against the previous $t\ll k_{\max}$ vectors only. For the first few experiments, the truncation parameter is fixed as~$t=2$. The sGMRES iteration terminates at the point when $\text{cond}(SAB_k)$ reaches or exceeds $10^{15}$. We refer to this method as ``sGMRES($10^{15}$)'' in the plots to follow.

\subsubsection{Convergence per outer cycle}

\begin{figure}
\centering\includegraphics[width=0.75\linewidth]{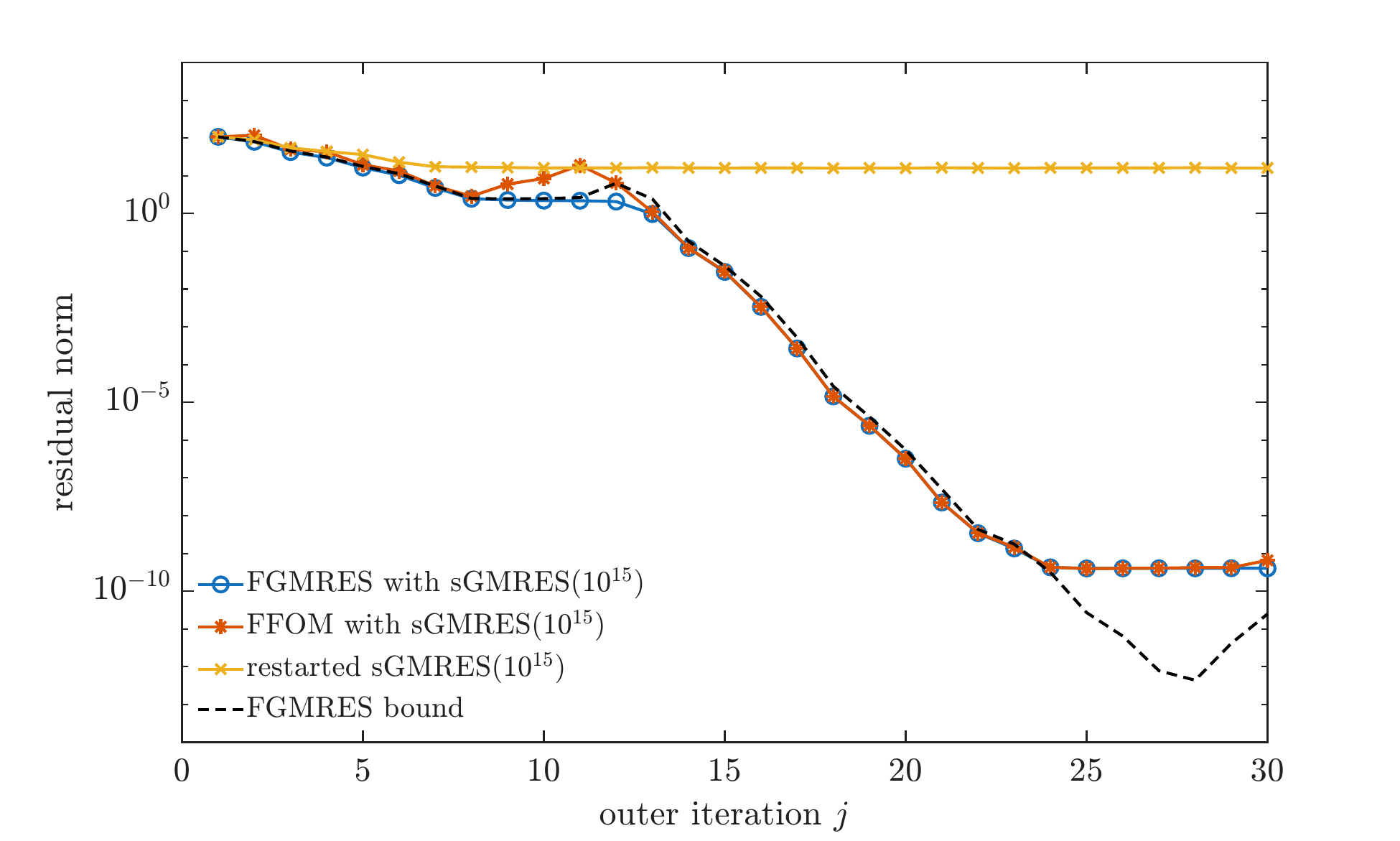}
\centering\includegraphics[width=0.75\linewidth]{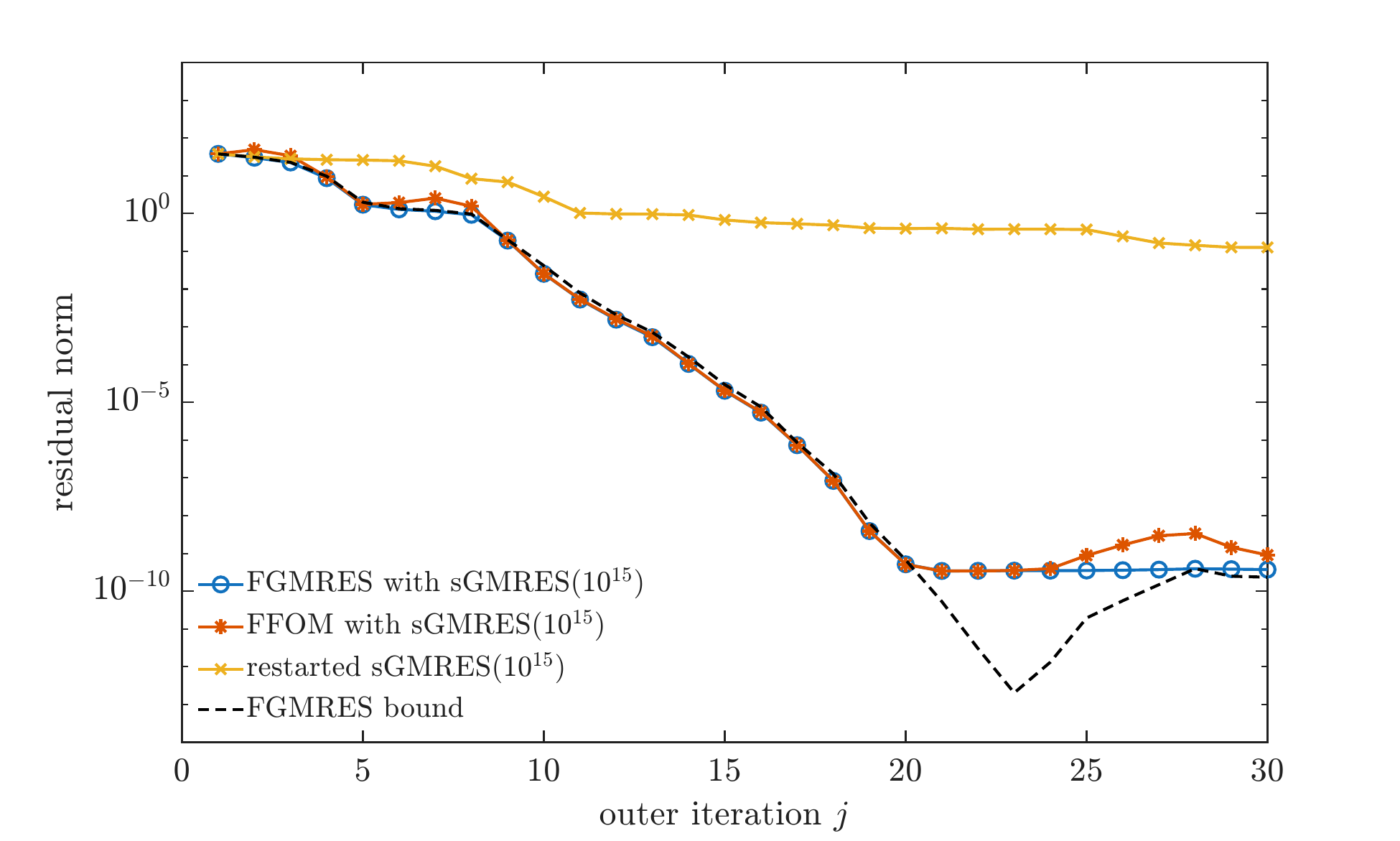}
    \caption{Residual convergence of FGMRES-stabilized sGMRES for the \texttt{vas\_stokes\_1M} problem. At each FGMRES iteration, the inner sGMRES solver is terminated when the sketched basis $S A B_k$ exceeds $10^{15}$ in condition number. We also show the convergence of restarted sGMRES using the same condition number criterion for restarting. The Arnoldi truncation parameter for sGMRES is $t=2$ (top) and $t=5$ (bottom). The dashed black line shows our upper bound~\eqref{eq:fgmres_bnd} on the FGMRES residual.}
    \label{fig:stokes_t2_res}
\end{figure}

In Figure~\ref{fig:stokes_t2_res} we compare FGMRES with sGMRES, FFOM with sGMRES, and restarted sGMRES in terms of the norm of the preconditioned residuals after $j$ outer iterations. 
The sGMRES bases $B_k$ generated during the $30$ FGMRES outer iterations vary in dimension between $115 \leq k\leq 135$. For restarted sGMRES, we observe $110 \leq k\leq 135$. Hence, over all iterations, the criterion $k\leq k_{\max}=500$ is satisfied, and so the inner sGMRES iterations are entirely controlled by the condition numbers of the generated Krylov bases, $\mathrm{cond}(SAB_k)$. Unfortunately, these condition numbers grow rather quickly, and so the values used for $k$ are too small for restarted sGMRES to make any significant progress.

Figure~\ref{fig:stokes_t2_res} also shows our FGMRES residual bound~\eqref{eq:fgmres_bnd}, which very closely tracks the actual FGMRES residual until FFOM stagnates. Upon stagnation of FFOM, the implicit formula for the FFOM residual~\eqref{eq:ffomres} is no longer reliable in floating point arithmetic and, in this case, underestimates the true residual. Overall, the residual bound from Section~\ref{sec:bnd} is well suited to terminate the inner FGMRES preconditioner (in this case sGMRES) as soon as the upper bound~\eqref{eq:fgmres_bnd} falls below the target residual tolerance. 

\subsubsection{Time-to-solution}

\begin{figure}
\centering\includegraphics[width=0.75\linewidth]{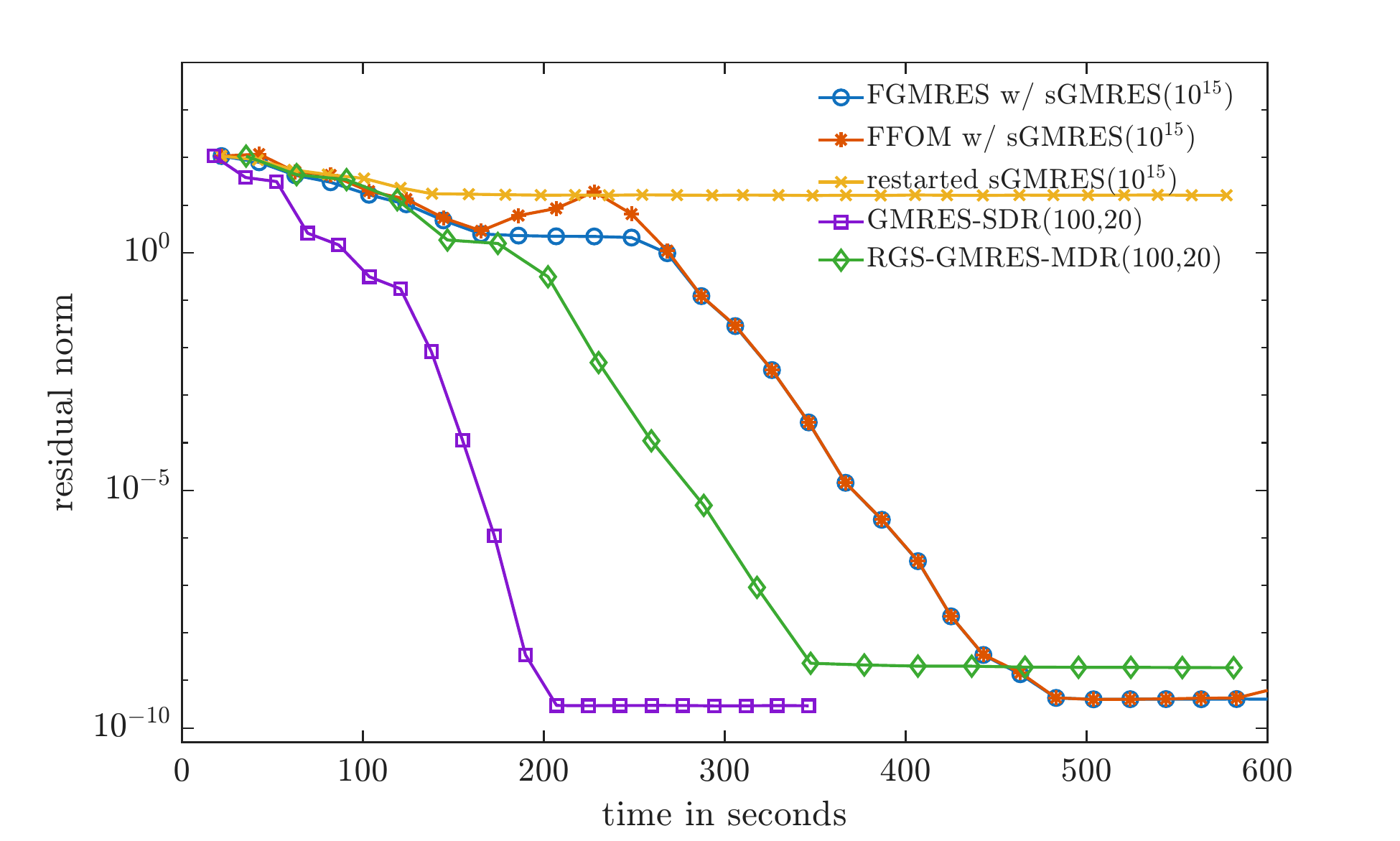}
\centering\includegraphics[width=0.75\linewidth]{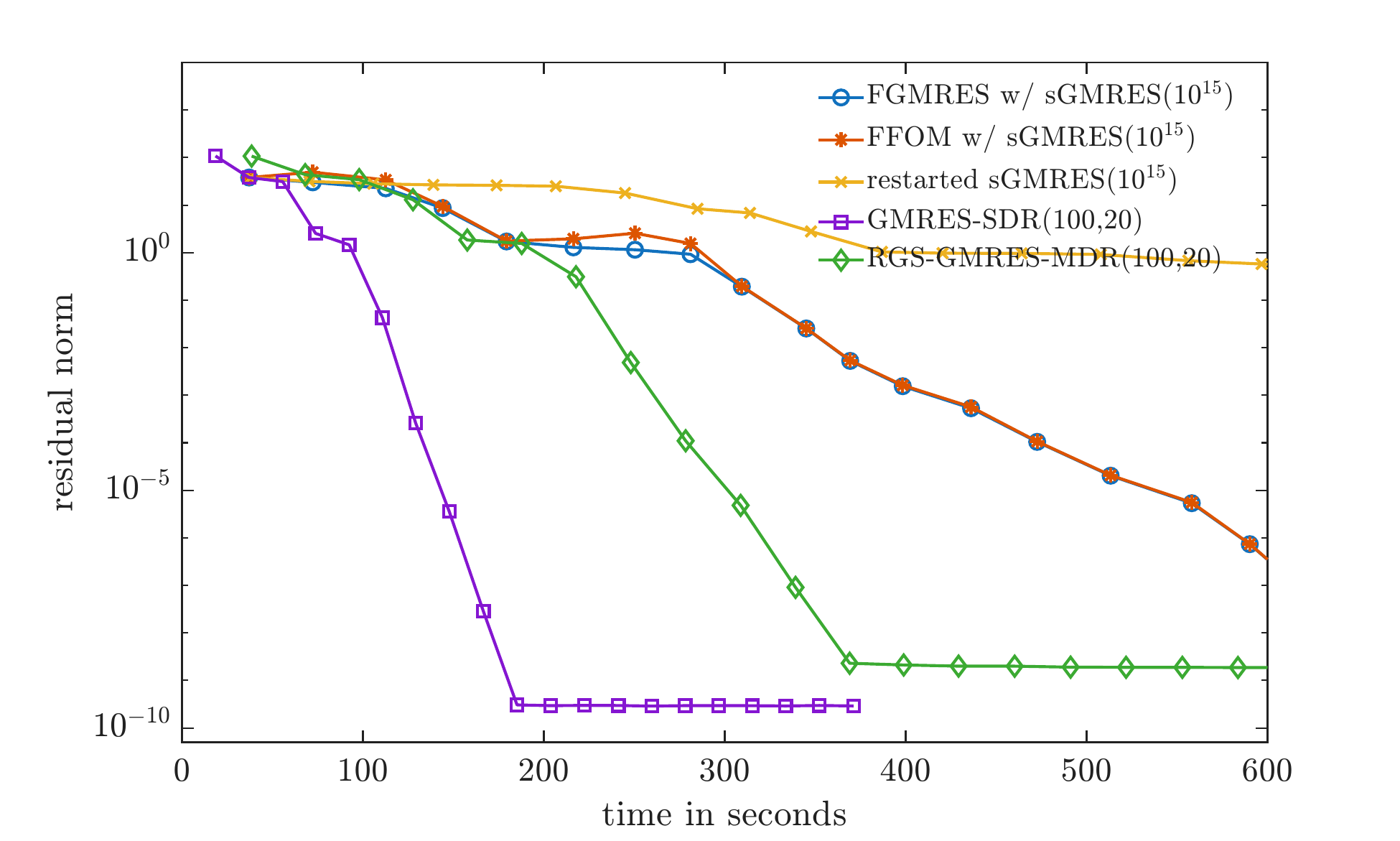}
    \caption{Wall-clock time versus residual of FGMRES-stabilized sGMRES for the \texttt{vas\_stokes\_1M} problem, compared to some other randomized methods. The Arnoldi truncation parameter for sGMRES and GMRES-SDR is $t=2$ (top) and $t=5$ (bottom).}
    \label{fig:stokes_t2_time}
\end{figure}

In terms of runtime per outer iteration, the FGMRES/FFOM and restarted variants are very similar, as shown in Figure~\ref{fig:stokes_t2_time}. The outer FGMRES/FFOM iteration does not add noticeable overhead, and both methods produce significantly smaller residuals than restarted sGMRES in less computation time.

In Figure~\ref{fig:stokes_t2_time} we compare with the GMRES-SDR (GMRES with sketching and deflated restarting) method from~\cite{BGS23}. The code is available online\footnote{\url{https://github.com/nla-group/GMRES-SDR} as of 2025/06/18}. We use the same parameters as in~\cite{BGS23}, namely a fixed Krylov basis dimension $k_\mathrm{fixed}=100$ for the inner iterations, and deflate up to $20$ harmonic Ritz vectors, which are improved from restart to restart.  All the other parameters (like sketching dimension~$s$ and Arnoldi truncation parameter~$t$) are chosen the same as for sGMRES above. We find that GMRES-SDR converges faster than FGMRES/FFOM on this problem (when $t=2$), but the latter methods have the advantage of being easier to implement and requiring no choice of $k_\mathrm{fixed}$ or the maximal number of harmonic Ritz vectors. 

In Figure~\ref{fig:stokes_t2_time} we also compare FGMRES-sGMRES with the \emph{randomized GMRES with SVD-based deflated restarting method}, RGS-GMRES-MDR($\overline{m}, \overline{k}$), one of the preferred methods in \cite{jang2025randomized}. Here, $\overline{m}$ refers to the maximal number of Krylov vectors to be stored per cycle and $\overline{k}$ is the number of approximate eigenpairs to be deflated. The code is available online\footnote{\url{https://github.com/Yongseok7717/RandomizedGMRES} as of 2025/06/18}.
We use the parameters $\overline{m}=100$ and $\overline{k}=20$, and also the same sketching operator (Clarkson--Woodruff with $s=1,000$) as in the other experiments. RGS-GMRES-MDR appears to make more progress per restart cycle than our FGMRES-sGMRES combination (although it stagnates on a slightly higher residual level), but overall both methods are comparable in terms of time-to-solution. 
We also highlight that for the FGMRES-sGMRES and FFOM-sGMRES approaches, there is a disadvantage in increasing the Arnoldi truncation parameter from $t=2$ to $t=5$ in terms of computation time, despite the marginal reduction in iteration numbers achieved (see Figure \ref{fig:stokes_t2_res}).

We highlight that RGS-GMRES-MDR, which is based on the ``RGS flexible Arnoldi iteration'' \cite[Alg.~1]{jang2025randomized}, uses \emph{two} applications of the sketching operator per inner iteration, while sGMRES can be implemented with just a single sketch per iteration. With the fast Clarkson--Woodruff sketch this is not very apparent in the total runtime, but more expensive sketches would make  RGS-GMRES-MDR significantly slower than FGMRES-sGMRES. We observed this in particular with the dense sketching operator used in the original implementation. On top of that, the FGMRES-sGMRES method requires no choice of $\overline{m}$ and $\overline{k}$ in advance.

\subsubsection{Influence of the Arnoldi truncation parameter}

\begin{figure}
    \centering\includegraphics[width=0.9\linewidth]{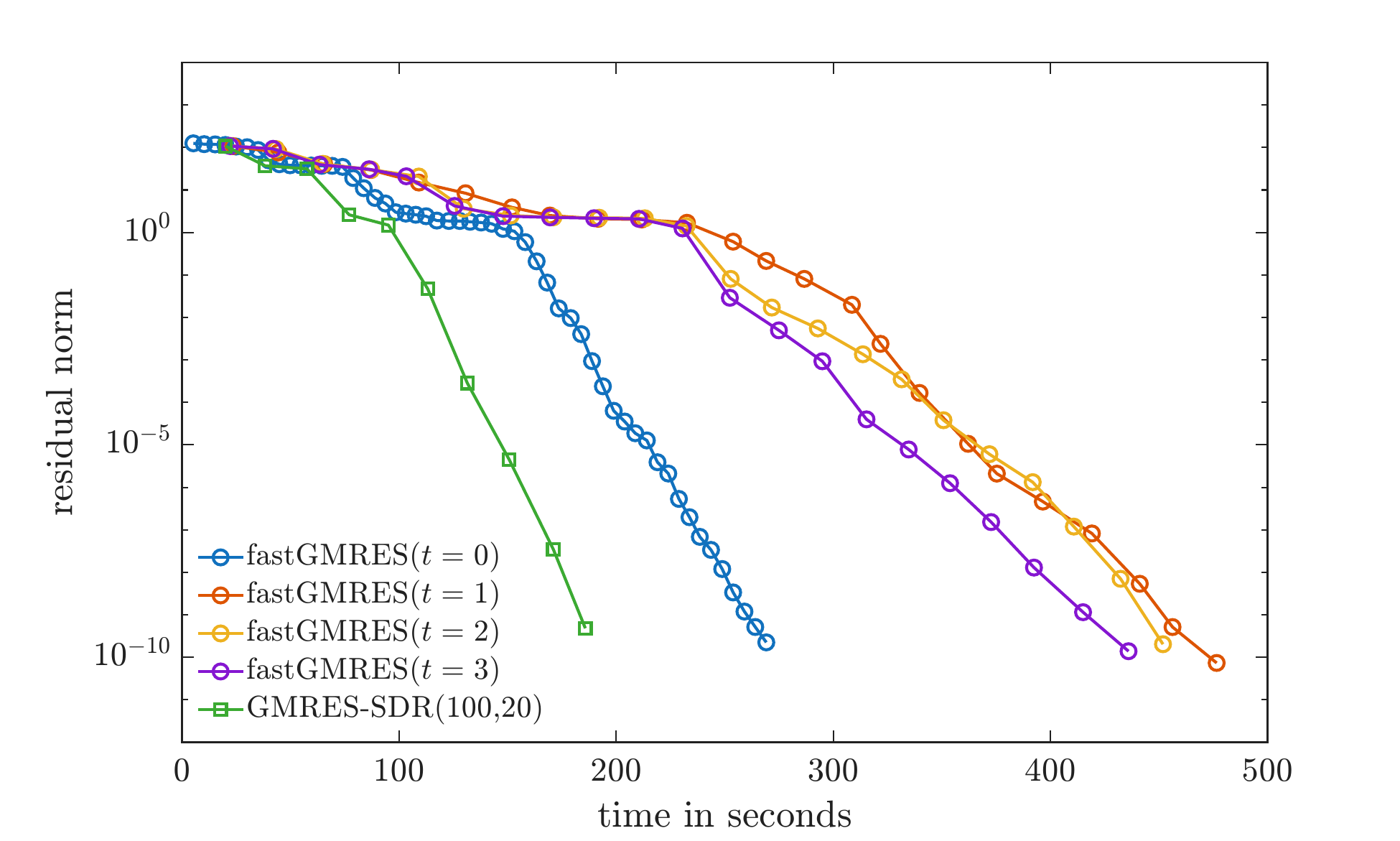}
    \caption{Exploring different values of the Arnoldi truncation parameter~$t$ for the \texttt{vas\_stokes\_1M} problem. The stabilizing effect of FGMRES can be clearly seen, with fastGMRES converging irrespective of the choice of~$t$. Surprisingly, $t=0$ (no orthogonalization at all) is the best choice for this problem, bringing fastGMRES close in performance to GMRES-SDR (using $t=2$).}
    \label{fig:stokes_k500_tx_time}
\end{figure}

In Figure~\ref{fig:stokes_k500_tx_time} we now explore the effect of the Arnoldi truncation parameter~$t$ on the runtime of our FGMRES-sGMRES combination. (Until now, this parameter has always been $t=2$.) For brevity, we now refer to our method as \emph{fastGMRES} in the plots (flexible and sketched GMRES with Arnoldi truncation), and we indicate the truncation parameter in parentheses. Generally, increasing $t$ can be expected to lead to better conditioned Krylov bases~$B_k$, meaning that the number of inner sGMRES iterations~$k$ per cycle increases. Decreasing $t$ will generally lead to worse conditioned bases, leading to a larger number of shorter sGMRES cycles. As can be seen in Figure~\ref{fig:stokes_k500_tx_time}, this intuition is confirmed by numerical evidence. It becomes obvious that FGMRES has a stabilizing effect on sGMRES: fastGMRES even converges for $t=0$, i.e., when $B_k$ is effectively a power basis without any orthogonalization! In this case, the number of Krylov basis vectors varies only between $28\leq k \leq 31$ per sGMRES cycle. 

We are surprised that $t=0$ is also the choice that leads to the best time-to-solution performance, coming closer to GMRES-SDR($100,20$). Whether choosing $t=0$ is a good idea will depend on the problem and implementation details. Either way, it is very encouraging that $t$ remains the only crucial parameter and fastGMRES performs robustly regardless of its choice.

\subsection{Problem \texttt{NACA12}}

We now consider a NACA airfoil problem, \texttt{NACA12},  with a nonsymmetric matrix~$A$ of size $n=81,920$ and given right-hand side vector~$b$; see \cite{jacobs1933characteristics}. 
This problem has been used in \cite{jang2025randomized} with ILU preconditioning to demonstrate the RGS-GMRES-MDR method.

\begin{figure}
\centering\includegraphics[width=0.52\linewidth]{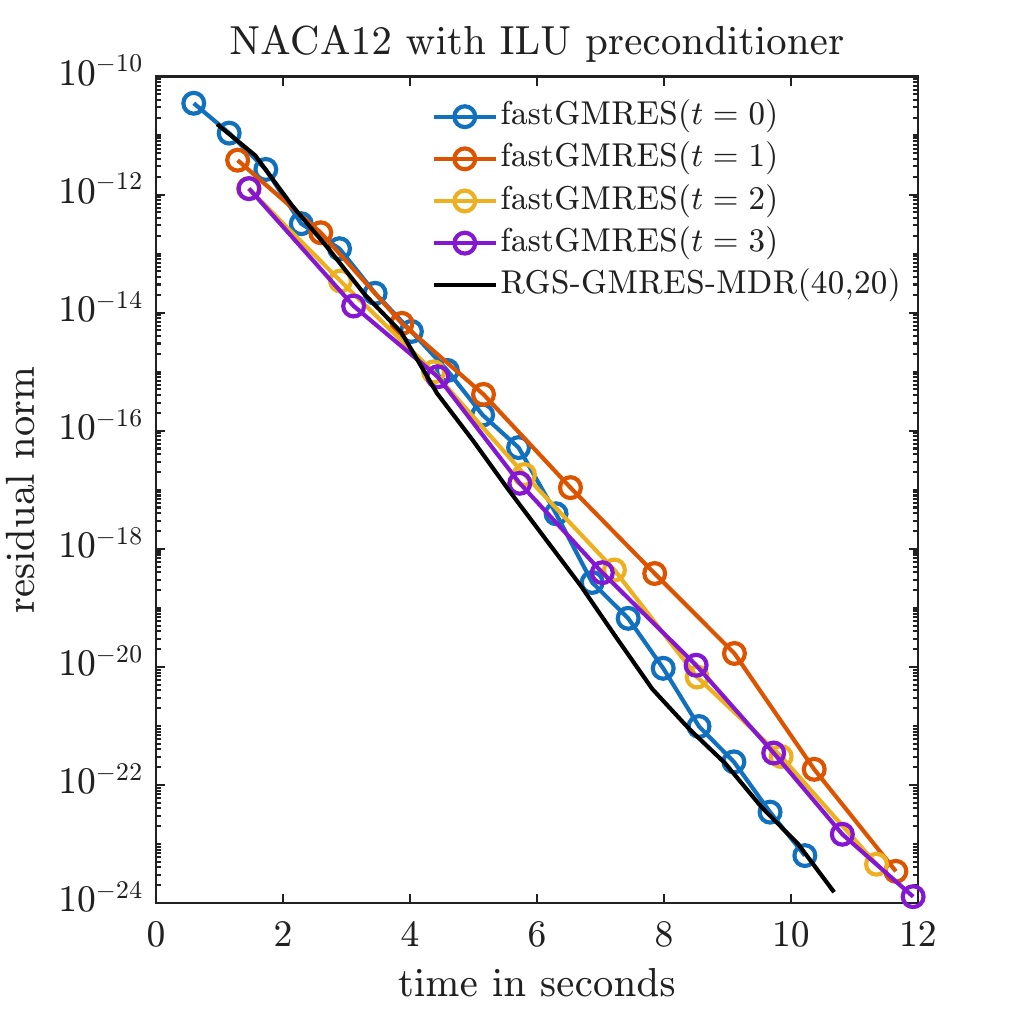}\hspace*{-2mm}\includegraphics[width=0.52\linewidth]{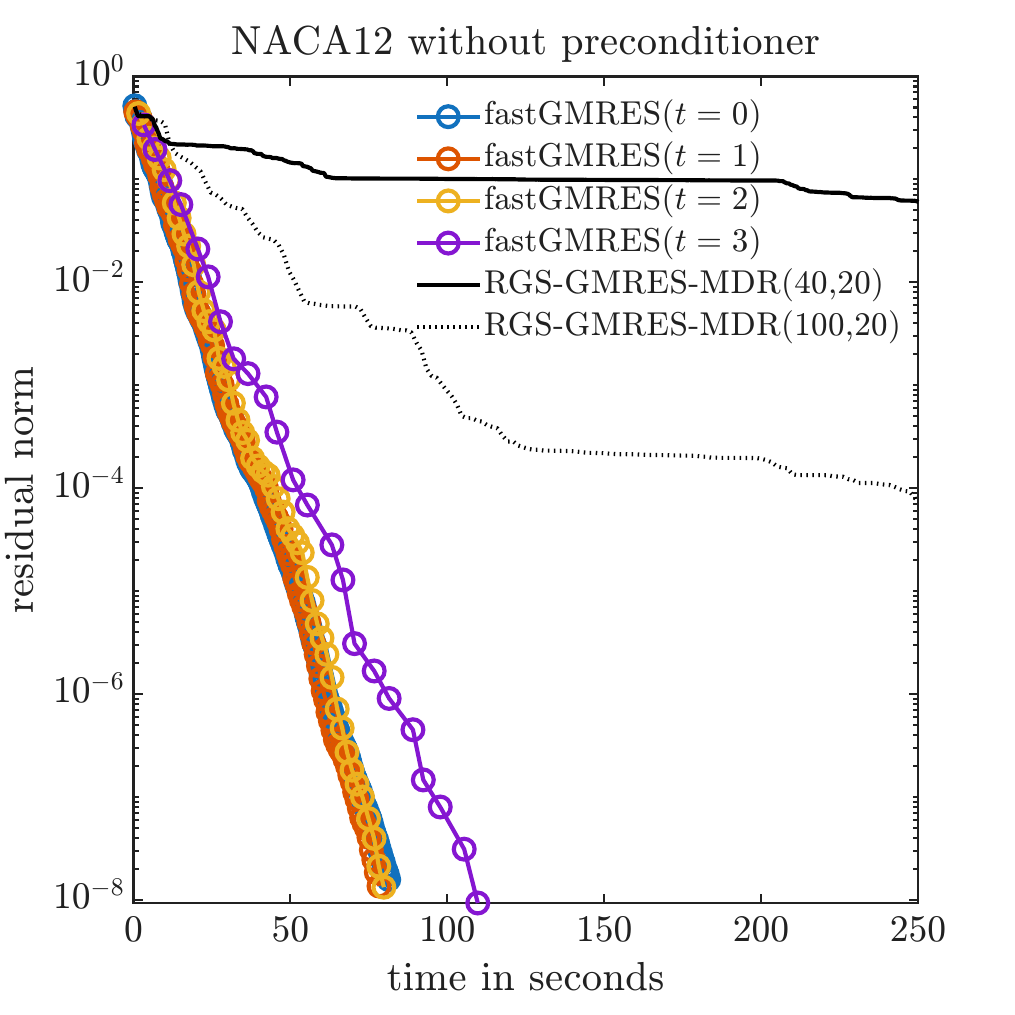}
    \caption{Exploring different values of the Arnoldi truncation parameter~$t$ for the \texttt{NACA12} problem, with (left) and without ILU preconditioning (right).}
    \label{fig:naca}
\end{figure}

We test the convergence of fastGMRES on the (zero-fill) ILU left-preconditioned problem, for different choices of the Arnoldi truncation parameter $t=0,1,2,3$. The results are shown on the left of  Figure~\ref{fig:naca}. We find that all tested choices for~$t$ work almost equally well, and overall this is an easy problem to solve due to the effectiveness of the ILU preconditioner. We also show the performance of RGS-GMRES-MDR($40,20$), which is essentially the same experiment as in \cite[Fig.~3]{jang2025randomized}, except that we have exchanged the dense sketching matrix for a sparse Clarkson--Woodruff sketch ($s=1,000$) to improve the performance. The RGS-GMRES-MDR method performs very similarly to fastGMRES on this preconditioned problem. 

On the right of Figure~\ref{fig:naca} we solve the same $Ax=b$ problem  but \emph{without} preconditioning. This becomes much harder but, again, fastGMRES with truncation parameter $t=0$ is among the fastest solvers for this problem. RGS-GMRES-MDR($40,20$) now converges very slowly, so we have also tried RGS-GMRES-MDR($100,20$) but timed it out after 5~minutes as the target relative residual norm of $10^{-8}$ was not reached.

\subsection{Problem \texttt{ML\_Geer}}

\begin{figure}
\centering\includegraphics[width=0.52\linewidth]{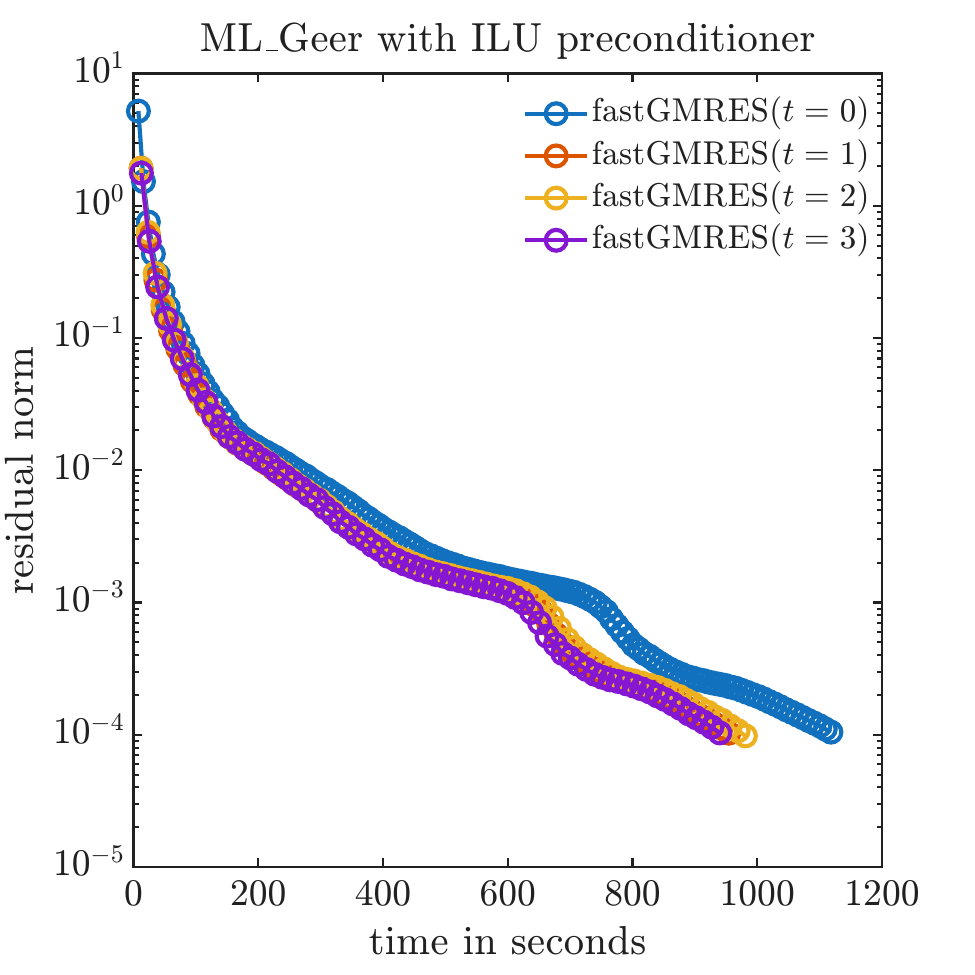}\hspace*{-2mm}\includegraphics[width=0.52\linewidth]{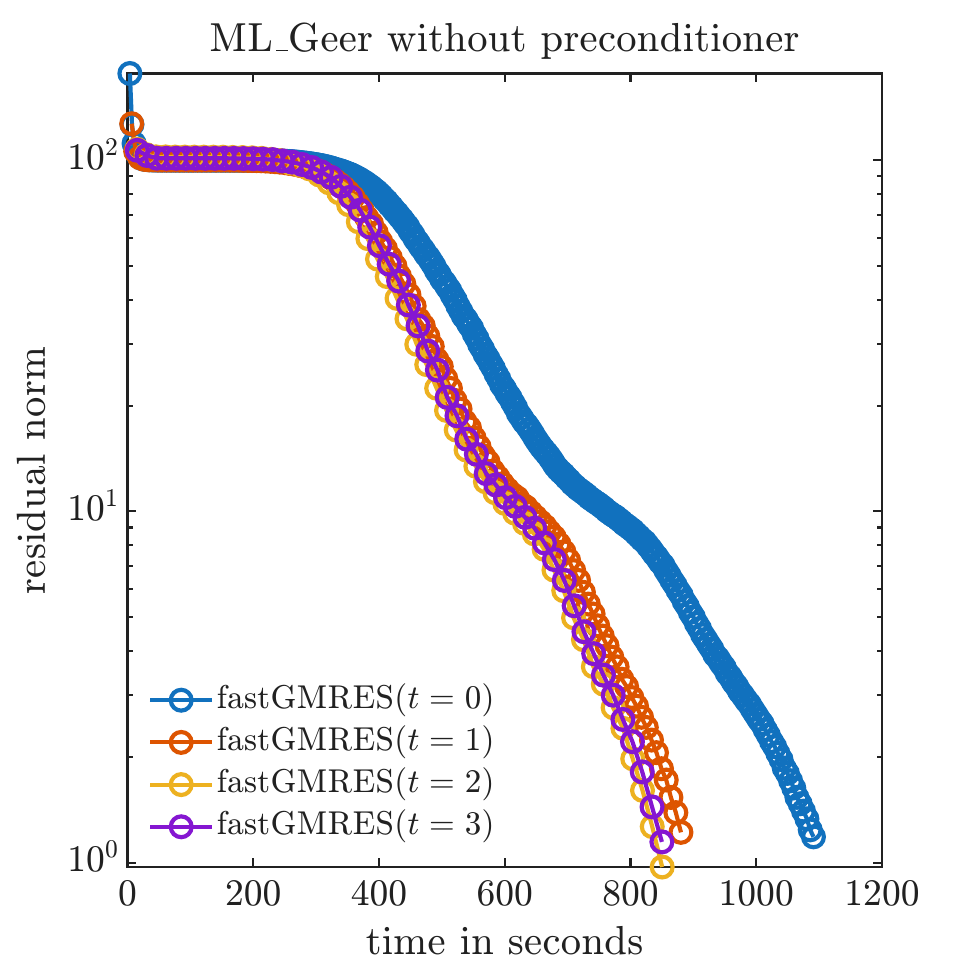}
    \caption{Exploring different values of the Arnoldi truncation parameter~$t$ for the \texttt{ML\_Geer} problem, with (left) and without ILU preconditioning (right).}
    \label{fig:ml_geer}
\end{figure}

We consider the \texttt{ML\_Geer} problem from the SuiteSparse Matrix Collection~\cite{davis2011university}, a nonsymmetric $A$ of size $n=1,504,002$, using a random right-hand side vector for $Ax=b$.
In Figure~\ref{fig:ml_geer} we solve that problem with and without (zero-fill) ILU left-preconditioner, with a target relative residual norm of $10^{-6}$ and $10^{-3}$, respectively, varying the Arnoldi truncation parameter~$t$. Again, we observe that the method converges robustly (albeit slowly) irrespective of the choice of $t$, but there is a noticeable delay in convergence with $t=0$ for later iterations. This delay is partly explained by the orthogonalization cost for the $m=137$, respectively $m=257$, outer FGMRES iterations that are performed. As the cost of full orthogonalization increases like $O(nm^2)$, this becomes noticeable for large~$m$ in particular when matrix--vector products are relatively cheap (which is the case when no preconditioner is used).

\subsection{Benchmarks against classical solvers}

\begingroup
\footnotesize
\tabcolsep = 1pt

\def\arraystretch{1.0}

\newcolumntype{C}{>{\centering\let\newline\\\arraybackslash\hspace{0pt}}m{1.3cm}}

\begin{table}
\caption{Runtimes and residuals (in parentheses) of various methods. The methods are fastGMRES, BiCGSTAB, GMRES (restart length 50 and 100), GMRES-DR (50 and 100 with 10 and 20 deflated vectors, respectively), and GMRES-SDR (50 and 100 with 10 and 20 deflated vectors, respectively). Each problem is run without and with preconditioning, the latter indicated by the suffix ``-p'' in the name. The target relative residual is $10^{-6}$ and methods are terminated when the runtime exceeds~300 seconds.\label{tab:bench}}
\begin{center}
\begin{footnotesize}
\begin{tabular}{lCCCCCCCC}\toprule
  \multirow{2}{*}{Problem}  &   \multirow{2}{1.3cm}{\centering fast\\ GMRES} &  \multirow{2}{1.3cm}{\centering BiCG\\ STAB} & \multicolumn{2}{c}{GMRES}  & \multicolumn{2}{c}{GMRES-DR}  & \multicolumn{2}{c}{GMRES-SDR} \\ 
    \cmidrule(lr){4-5}\cmidrule(lr){6-7}\cmidrule(lr){8-9}
   &  & & (50) & (100) & (50,10) & (100,20) & (50,10) & (100,20) \\[1mm]\midrule\\[-3mm]

atmosmodd & 10.7 (5.0e-7) & \textbf{8.1} (6.5e-7) & 21.5 (1.0e-6) & 30.6 (1.0e-6) & 22.6 (9.6e-7) & 29.6 (9.9e-7) & 14.1 (7.4e-7) & 13.6 (7.5e-7)  \\[4mm]

atmosmodd-p & 9.8 (1.1e-7) & 8.8 (6.7e-7) & 11.2 (9.7e-7) & 12.5 (9.9e-7) & 10.0 (8.3e-7) & 10.0 (9.9e-7) & 8.6 (7.0e-7) & \textbf{8.4} (8.1e-7)  \\[4mm] 

NACA12 & 54.3 (9.5e-7) & \cancel{7.2} (4.9e-1) & 300.0 (2.5e-1) & 300.0 (1.3e-1) & 300.0 (4.6e-3) & 79.6 (1.0e-6) & 300.4 (4.8e-3) & \textbf{53.2} (8.2e-7)  \\[4mm]

NACA12-p & 5.5 (9.4e-7) & 7.2 (7.3e-7) & 7.6 (9.5e-7) & 5.8 (9.6e-7) & 4.3 (9.3e-7) & 4.2 (9.9e-7) & 4.3 (6.8e-7) & \textbf{3.9} (6.7e-7)  \\[4mm]

SiO2 & 27.1 (6.0e-7) & 36.9 (8.5e-7) & 167.2 (1.0e-6) & 78.7 (9.9e-7) & 10.3 (1.0e-6) & 10.8 (1.0e-6) & 12.2 (7.8e-7) & \textbf{8.5} (7.6e-7)  \\[4mm] 

SiO2-p & 20.2 (8.0e-7) & 16.9 (6.5e-7) & 9.6 (9.5e-7) & 6.2 (9.4e-7) & 6.7 (9.3e-7) & \textbf{5.6} (9.8e-7) & 12.2 (6.6e-7) & 18.2 (3.0e-7)  \\[4mm] 

t2em & 41.0 (9.0e-7) & \textbf{29.5} (9.8e-7) & 300.1 (1.1e-5) & 300.1 (8.7e-6) & 115.3 (1.1e-6) & 186.4 (1.7e-6) & 62.7 (8.6e-7) & 60.6 (9.1e-7)  \\[4mm]

t2em-p & 32.2 (8.7e-7) & \textbf{29.8} (9.5e-7) & 92.0 (1.0e-6) & 75.7 (9.9e-7) & 45.6 (9.9e-7) & 48.6 (1.0e-6) & 34.3 (7.1e-7) & 33.0 (7.6e-7)  \\[4mm]

vas\_stokes\_1M & \textbf{301.8} (4.6e-5) & \cancel{129.0} (4.6e-2) & 300.2 (2.3e-2) & 300.1 (2.2e-2) & 300.5 (5.1e-3) & 300.6 (4.6e-3) & 301.2 (1.8e-2) & 300.5 (2.2e-2)  \\[4mm]

vas\_stokes\_1M-p & 168.7 (8.1e-7) & 300.3 (1.1e-3) & 300.2 (1.9e-2) & 300.3 (7.1e-3) & 165.3 (9.9e-7) & 134.7 (9.5e-7) & 153.0 (9.2e-7) & \textbf{124.9} (7.3e-7)  \\[4mm]

\bottomrule
\end{tabular}
\end{footnotesize}
\end{center}
\end{table}

\endgroup

We now compare fastGMRES to a number of other solvers, including restarted GMRES~\cite{SaadSchultz1986},  BiCGSTAB~\cite{VanDerVorst1992}, GMRES with deflated restarting (GMRES-DR)~\cite{Morgan2002}, and restarted GMRES with sketching and deflated restarting (GMRES-SDR)~\cite{BGS23}. Our approach is as follows. For every problem we consider in this test, we solve $Ax=b$ either with or without (zero-fill) ILU right-preconditioning. In all cases, the target relative residual norm is $10^{-6}$ and all methods are timed out after 300 seconds. 

In Table~\ref{tab:bench}, we report the achieved residual norm as well as the total runtime. The results are mixed, but there are a number of general observations we can make. First, on 3 of the 10 problems, BiCGSTAB is the fastest solver, but it is also the least reliable as it fails on 2/10 problems (with \texttt{flag=4} in MATLAB, which means that one of the scalar quantities calculated by the BiCGSTAB algorithm became too small or too large to continue computing). 

GMRES-SDR(100,20) also performs very well, being the fastest solver for 5/10 problems. Generally, deflated restarting helps a lot in mitigating convergence delays caused by restarting. Even without sketching, GMRES-DR shows good performance, but for all deflated restarting methods the selection of parameters (subspace size and number of deflated eigenvectors) is crucial for success and performance. For example, for \texttt{NACA12}, GMRES-DR(50,10) and GMRES-SDR(50,10) fail to converge within 300 seconds, while GMRES-DR(100,20) and GMRES-SDR(100,20) converge in 79.6 seconds and 53.2 seconds, respectively. While the sketching in GMRES-SDR typically yields a performance gain over GMRES-DR, sketching introduces additional parameters such as the Arnoldi truncation parameter $t$ and the embedding dimension~$s$ (here, these are chosen as $t=2$ and $s=1,000$). 

Our new method, fastGMRES, provides an attractive alternative in that it is among the most reliable methods, only failing to converge within 300 seconds on one of the problems (\texttt{vas\_stokes\_1M}) for which all other methods are timed out as well. Among all solvers applied to this problem, fastGMRES achieves the smallest relative residual (4.6e-5) within 300 seconds. At the same time, for our experiments fastGMRES requires essentially no parameter choices. In all cases, we have kept $t=0$ and $s=1,000$ fixed and fastGMRES completes all problems robustly, in the case of \texttt{vas\_stokes\_1M} beyond the 300 second limit. While fastGMRES does not outperform the fastest solver (such as GMRES-SDR with carefully chosen parameters) on any of the other problems, it is often close in performance and reliable (as opposed to BiCGSTAB). For 9/10 problems, fastGMRES outperforms  GMRES(50) and GMRES(100).

\section{Conclusions and future work}\label{sec:concl}

We have demonstrated that there is benefit in using sGMRES as a preconditioner within FGMRES, resulting in a very simple but fast randomized solver that is guaranteed to provide non-increasing residual norms in exact arithmetic. The main parameter to choose is the Arnoldi truncation parameter~$t$, and we have been surprised that a competitive choice in most experiments is $t=0$; that is, to generate a power basis without any orthogonalization. This means that no inner products are being computed within sGMRES. In future work, we would like to perform a more comprehensive performance comparison on an even wider range of problems, benchmarking both deterministic and randomized Krylov solvers. We are currently developing a benchmark collection that focuses on problems from PDEs and optimization for that purpose.

While more work is needed to fully understand and justify why using $t=0$ as the truncation parameter in our method leads to fast serial runtime, the resulting method may also be advantageous in view of other computational considerations. Most notably, as is also mentioned, e.g., in \cite{brown2025inner}, avoiding inner products as much as possible can be beneficial when solving problems with mixed precision and for parallel computing. Requiring less fine-grained communication for orthogonalization means, for example, that one may consider precomputing a few matrix powers $A,A^2,\ldots\,$ explicitly and applying these in parallel to compute the inner sGMRES basis (and its sketches) faster. If $A$ is sufficiently sparse with matrix powers that do not fill-in too quickly, this may indeed be a viable route to higher parallelism.

Finally, we recall that for the \texttt{ML\_Geer} problem in Figure~\ref{fig:ml_geer}, the number of outer FGMRES iterations~$m$ gets rather large (for $t=0$, $m$ increases to $137$ and $257$, respectively). In this regime, the increasing orthogonalization cost for the FGMRES basis~$V_m$ becomes noticeable. One obvious solution would be to restart the outer FGMRES loop. However, as we have seen in previous experiments (see, e.g., Figure~\ref{fig:stokes_t2_res}), restarting may lead to stagnation or significant delay in convergence.
It might therefore be desirable not to restart, but instead to reduce the arithmetic cost of the FGMRES iteration. One possible approach is to run the outer FGMRES iteration only up to a limited basis size $m_\mathrm{max}$, and then to interpret the resulting method as the inverse of a preconditioner $M_j$ for a higher-level FGMRES loop. This process may be recursively applied, leading to a ``multi-level FGMRES method''. We hope to explore these ideas in future work.

\section*{Acknowledgments}
Both authors are grateful to two anonymous referees for their valuable comments. We also thank Lauri Nyman for useful input. 
Both authors acknowledge joint funding from the UK's Engineering and Physical Sciences Research Council (EPSRC grant EP/Z533786/1). SG is supported by the Royal Society (RS Industry Fellowship IF/R1/231032).

\bibliographystyle{plain}
\bibliography{refs}

\end{document}